\documentclass[12pt]{amsart}
\usepackage{amscd,mathpazo,fullpage}
\usepackage{graphicx} 
\usepackage{amsmath}
\usepackage{amsfonts}
\usepackage{amssymb}
\usepackage{amsmath} 
\usepackage{stmaryrd}

\usepackage{geometry}
\usepackage{mathabx}
\usepackage{tikz-cd}
\usepackage{xcolor}
\usepackage{array}

\newcommand{\lb}{\left(}
\newcommand{\rb}{\right)}

\newcommand{\cg}{\nabla^{g}}
\newcommand{\ch}{\nabla^{h}}
\newcommand{\ck}{\nabla^{K}}
\newcommand{\cc}{\nabla^{c}}
\newcommand{\cb}{\nabla^{B}}

\newcommand{\X}{\mathfrak{X}}
\newcommand{\id}{\operatorname{Id}}
\newcommand{\dd}{\mathrm{d}}

\newtheorem{theorem}{Theorem}
\newtheorem{corollary}[theorem]{Corollary}
\newtheorem{lemma}[theorem]{Lemma}
\newtheorem{proposition}[theorem]{Proposition}

\theoremstyle{definition}
\newtheorem{definition}[theorem]{Definition}
\newtheorem{example}[theorem]{Example}
\newtheorem{remark}[theorem]{Remark}

\title[Born geometry via K\"unneth structures and recursion operators]{Born geometry via K\"unneth structures and\\ recursion operators}
\author{M.~J.~D.~Hamilton}
\address{Mathematisches Institut, {\smaller LMU} M\"unchen, Theresienstr.~39, 80333~M\"unchen, Germany}
\email{mark.hamilton@math.lmu.de}
\author{D.~Kotschick}
\address{Mathematisches Institut, {\smaller LMU} M\"unchen, Theresienstr.~39, 80333~M\"unchen, Germany}
\email{dieter@math.lmu.de}
\author{P.~N.~Pilatus}
\address{Fachbereich Mathematik, Universit\"at Hamburg, Bundesstr.~55, 20146 Hamburg, Germany}
\email{paula.pilatus@uni-hamburg.de}

\date{\today; \copyright{\ M.~J.~D.~Hamilton, D.~Kotschick and P.~N.~Pilatus 2024}}

\begin{document}

\begin{abstract}
We propose a simple definition of a Born geometry in the framework of K\"unneth geometry. While superficially 
different, this new definition is equivalent to 
the known definitions in terms of para-quaternionic or generalized geometries. We discuss integrability of Born structures
and their associated connections. In particular we find that for integrable Born geometries the Born connection is obtained
by a simple averaging under a conjugation from the K\"unneth connection. We also give examples of integrable Born geometries on nilmanifolds.
\end{abstract}

\maketitle


\section{Introduction}

Born geometry was introduced by Freidel, Leigh and Minic in~\cite{freidel2014born} in the context of T-duality in string theory. These authors 
proposed to describe T-duality by replacing the target space of a theory by a manifold of twice the dimension, carrying a geometric structure 
which they call a Born structure or a Born geometry. The name was chosen because Max Born~\cite{born1938suggestion} had suggested a long time ago that, 
in order to unify quantum mechanics and general relativity, the momentum space should be allowed to have curvature too. Born geometry
and its application in high-energy physics have been developed further in many papers, see e.g.~\cite{freidel2014quantum,freidel2015metastring,freidel2017generalised,freidel2019unique,marotta2018parahermitian,marotta2021born}.

The mathematical definition of a Born geometry that appears for example in the paper of Freidel, Rudolph and Svoboda~\cite{freidel2019unique} 
was in terms of para-quaternionic and 
para-hermitian structures, and was influenced by the setup of generalized geometry. An equivalent description of Born structures which is 
purely in the language of generalized geometry is discussed in~\cite{hu2019zro}.

In this paper we propose a streamlined definition of a Born geometry as a diagram 
    \begin{center}
        \begin{tikzcd}[column sep=small]& g\arrow[dl,"A"']\arrow{dr}{B} & \\\omega \arrow{rr}{-J}   & &h \ \ ,
     \end{tikzcd}
    \end{center}
    where $g$ and $h$ are pseudo-Riemannian metrics and $\omega$ is a non-degenerate 2-form, such that the 
    intertwining operators satisfy $A^2=B^2=-J^2=\operatorname{Id}$. We do not assume that the operators define any kind
    of (para-)quaternionic structure, since this turns out to be automatically true. Further, we make no assumption about the 
    signatures of the pseudo-Riemannian metrics $g$ and $h$. It turns out that $g$ is automatically of neutral signature,
    but the signature of $h$ is arbitrary as long as it is of the form $(2p,2q)$.
    
This definition of a Born structure is motivated by the definition of a hypersymplectic structure in~\cite{bande2008geometry}.
This structure was originally defined by Hitchin~\cite{hitchin1990hypersymplectic} using a neutral metric and certain 
symplectic forms that are parallel for its Levi-Civita connection. It was shown in~\cite{bande2008geometry}, see 
also~\cite[Chapter~8]{kuennethgeom}, that Hitchin's definition is equivalent to giving only a diagram of symplectic forms
\begin{center}
        \begin{tikzcd}[column sep=small]& \omega\arrow[dl,"A"']\arrow{dr}{B} & \\\alpha \arrow{rr}{J} & & \beta
     \end{tikzcd}
    \end{center}
such that the intertwining operators satisfy $A^2=B^2= -J^2= \id$. All other relations between the operators
and the existence of a neutral metric defined in terms of these symplectic forms are a consequence of the diagram.

In our definition of a Born geometry it turns out that the $(\pm 1)$-eigenbundles of the involution $A$ are Lagrangian
for the form $\omega$, and so constitute an almost bi-Lagrangian~\cite{GS} or almost para-K\"ahler~\cite{Cr}
or almost K\"unneth structure~\cite{kuennethgeom}. 
As usual, the use of the almost qualification means that no integrability has been assumed. Conversely, we will see 
that every Born geometry is just an almost K\"unneth structure together with a choice of an isomorphism between
the two Lagrangian subbundles.

The existence of such almost structures is purely a matter of algebraic topology and bundle theory. 
However, once one imposes integrability, it turns out that these structures are geometrically interesting,
and their existence can be quite subtle. We will discuss the (partial) integrability of Born structures, and 
exhibit some consequences. For example, every almost K\"unneth structure on a manifold gives rise to 
a preferred affine connection, the K\"unneth connection, with respect to which the whole structure is parallel.
The K\"unneth connection is torsion-free if and only if the structure is integrable. Now Freidel, Rudolph and Svoboda~\cite{freidel2019unique} proved that every Born structure gives rise to a unique connection 
for which the structure is parallel, and which satisfies a variation of torsion-freeness suggested by generalized geometry.
The expressions given for the Born connection in~\cite{freidel2019unique} are very complicated, and are
difficult to work with. We will prove that in the integrable case the Born connection arises in a very simple 
way from the K\"unneth connection of the underlying K\"unneth structure, although in general it does not equal the 
K\"unneth connection.

\subsection*{Structure of the paper}

In Section~\ref{s:prelim} we recall some known facts, mostly about K\"unneth geometry, which we 
need for our discussion of Born geometry. In Section~\ref{s:Born} we discuss our definition of 
Born geometry, and we develop simple properties. We also discuss integrability and the relationship to 
hypersymplectic geometry. In Section~\ref{sect: connections} we compare the different connections 
that arise from a Born geometry, and we prove that in the integrable case the Born 
connection can be recovered directly from the K\"unneth connection. 
Finally in Section~\ref{s:examples}
we give examples of integrable Born geometries on some nilmanifolds.

\section{Preliminaries}\label{s:prelim}

In this section we recall some definitions we need for our formulation of Born geometry.
All the missing details are contained in~\cite{kuennethgeom}.

\subsection{K\"unneth structures}
\label{sect: def almost kunneth}

The following is our basic structure.
\begin{definition}\label{d:kunneth}
An \emph{almost K\"unneth structure} on a smooth manifold $M$ consists of a non-degenerate 2-form $\omega$ 
together with two complementary $\omega$-isotropic subbundles $F$ and $G$. An almost K\"unneth structure is
\emph{K\"unneth} if it is integrable, meaning that if $\omega$ is closed and therefore symplectic, and $F$ and $G$ 
are integrable to Lagrangian foliations $\mathcal{F}$ and $\mathcal{G}$.
\end{definition}
These structures are also known under several other names, the most common of which is that of a 
para-K\"ahler structure, cf.~\cite{Cr}. It is proved in~\cite[Section~6.3]{kuennethgeom} that an almost K\"unneth structure
is equivalent to an almost para-Hermitian structure. Moreover, in the integrable case, a K\"unneth structure 
is equivalent to a para-K\"ahler structure.

One can think of Definition~\ref{d:kunneth} as a metric-free definition of para-K\"ahler structures.
The metric is recovered by the following construction.
The splitting of the tangent bundle $TM=F\oplus G$ defines an almost product structure
\begin{align*}
    I\colon TM &\longrightarrow TM\\
    X_F+X_G&\longmapsto X_F-X_G \ ,
\end{align*}
which, together with $\omega$, yields:
\begin{proposition}\cite[Proposition~5.9]{kuennethgeom}
\label{prop: neutral metric}
For an almost K\"unneth structure $(\omega,F,G)$ the map 
\begin{align*}
    g\colon TM\times TM &\longrightarrow\mathbb{R}\\
    (X,Y)&\longmapsto \omega(IX,Y)
\end{align*}
defines a pseudo-Riemannian metric of neutral signature.
\end{proposition}

Note that, while in the definition of an almost K\"unneth structure the two subbundles are treated symmetrically, after defining 
the almost product structure, this symmetry is broken, and the subbundles are distinguished by the property of being either 
the $(+1)$- or the 
$(-1)$-eigenbundles of $I$. Therefore, the almost product structure and the neutral metric associated with an 
almost K\"unneth structure are well-defined only up to an overall sign.

There is a canonical connection associated with an almost K\"unneth structure, the \emph{K\"unneth connection},
whose definition goes back to Hess~\cite{H1,H2}.
\begin{theorem}\cite[Theorem~6.6]{kuennethgeom}
Let $(\omega, F,G)$ be an almost K\"unneth structure on a smooth manifold $M$. There exists a unique affine connection $\ck$ on $M$, the \emph{K\"unneth connection}, that preserves both $F$ and $G$, is compatible with $\omega$, and whose mixed torsion vanishes identically, i.e. $T^K(X,Y)=0$ for $X\in F$ and $Y\in G$ or vice versa, where $T^K$ is the torsion tensor of $\ck$.
\end{theorem}
The K\"unneth connection can be expressed as follows. Let 
$$
D\colon \X(M)\times\X(M)\longrightarrow\X(M)
$$ 
be the map defined by $i_{D(X,Y)}\omega =L_X i_Y\omega$.
Then for $X\in\X(M)$ and $Y\in\Gamma(F)$ the formula 
\begin{align*}
    \ck_X Y =D(X_F,Y)_F + [X_G,Y]_F
\end{align*}
defines a connection on $F$. Similarly, for $X\in \X(M)$ and $Y\in\Gamma(G)$
\begin{align*}
    \ck_X Y = D(X_G,Y)_G+[X_F,Y]_G
\end{align*}
defines a connection on $G$. Here, for $X\in\X(M)$, $X_F$ and $X_G$ denote the projections to $F$ and $G$, respectively. The K\"unneth connection is then given by
\begin{align*}
    \ck_X Y = \ck_X Y_F + \ck_X Y_G
\end{align*}
for $X,Y\in\X(M)$. In other words, it is the direct sum of the connections defined separately on 
$F$ and $G$.

The torsion of the K\"unneth connection is related to the integrability of the almost K\"unneth structure:

\begin{theorem}\cite[Theorem~6.8]{kuennethgeom}
An almost K\"unneth structure is K\"unneth if and only if its K\"unneth connection is torsion-free.
\end{theorem}
Since the K\"unneth connection is compatible with $\omega$ and commutes with the almost product structure $I$ defined by the 
subbundles $F$ and $G$, it is also compatible with the corresponding pseudo-Riemannian metric $g$. Therefore, we obtain the following:
\begin{theorem}\cite[Theorem~6.10]{kuennethgeom}
\label{thm: ck is LC}
If $(\omega, F,G)$ is a K\"unneth structure, then the Levi-Civita connection of the associated neutral metric $g$ 
is the K\"unneth connection.
\end{theorem}

\subsection{Recursion operators}

The almost product structure $I$ is a particular instance of what we call a \emph{recursion operator}. 
These will play a central role in our definition of Born structures. The terminology of recursion operators originates from 
the theory of bi-Hamiltonian systems, and was developed in a way that is relevant to this paper in~\cite{bande2008geometry,kuennethgeom}.

Given two non-degenerate bilinear forms $a$ and $b$ on the tangent bundle of a smooth manifold $M$, there is a unique field of invertible endomorphisms $A$, such that $a(A\cdot,\cdot)=b(\cdot,\cdot)$. Then, $A$ is called the \emph{recursion operator} from $a$ to $b$ and we will depict the situation as $a\overset{A}{\longrightarrow} b$. In this case, $A^{-1}$ is the recursion operator from $b$ to $a$, i.e. $b\overset{A^{-1}}{\longrightarrow} a$. Given a third bilinear form $c$, such that $b\overset{B}{\longrightarrow} c$, it follows that $a\overset{AB}{\longrightarrow} c$.

In the special case where a recursion operator satisfies $A^2=\operatorname{Id}$ and $A\neq \pm\operatorname{Id}$, the tangent
bundle $TM$ splits as the direct sum of the $(\pm 1)$-eigenbundles of $A$, see the discussion in Subsection~3.1 of~\cite{bande2008geometry}.

\section{Born geometry}\label{s:Born}

We can now formulate our definition of Born geometry.

\begin{definition}
\label{def: Born structure}
Let $M$ be a smooth manifold. A \emph{Born structure} on $M$ is a triple $(g,h,\omega)$, where $g$ and $h$ are pseudo-Riemannian metrics and $\omega$ is a non-degenerate 2-form, such that the recursion operators specified in the diagram
    \begin{center}
        \begin{tikzcd}[column sep=small]& g\arrow[dl,"A"']\arrow{dr}{B} & \\\omega \arrow{rr}{-J}   & &h
     \end{tikzcd}
    \end{center}
satisfy $A^2=B^2=-J^2=\operatorname{Id}$.
\end{definition}
Note that the existence of a non-degenerate 2-form implies that $M$ must be even-dimensional. 

A diagram of recursion operators as in this definition always commutes, because of the uniqueness of 
recursion operators. For example, $A\circ (-J)$ is a recursion operator from $g$ to $h$, and so is $B$.
Therefore, $-AJ=B$.

\subsection{Algebraic identities}

We will now discuss basic properties of Born structures that can be deduced immediately from the definition. 
In the following, let $M$ be a smooth manifold of dimension $2n$ with a Born structure $(g,h,\omega)$ with recursion operators 
$A$, $B$ and  $J$ as in Definition~\ref{def: Born structure}. 

\begin{lemma}
\label{lem: para-quaternionic}
The recursion operators $A$, $B$ and $J$ pairwise anti-commute and $ABJ=\operatorname{Id}$.
\end{lemma}
\begin{proof}
The second assertion follows immediately from $AB=-J$ and $J^2=-\id$. Moreover,
\begin{equation*}
AB=-J=J^{-1}=-B^{-1}A^{-1}=-BA \ .
\end{equation*}
This implies
\begin{equation*}
    AJ=-A^2B=ABA=-JA
\end{equation*}
and
\begin{equation*}
    JB=-AB^2=BAB=-BJ \ .
\end{equation*}
\end{proof}

An immediate consequence of Lemma \ref{lem: para-quaternionic} is the following:
\begin{lemma}
\label{lem: J_interchanges_esp}
The almost complex structure $J$ interchanges the $(\pm 1)$-eigenbundles of $A$. 
Similarly, it interchanges the $(\pm 1)$-eigenbundles of $B$. Furthermore, $A$ interchanges the eigenbundles of $B$ and vice versa.
\end{lemma}
\begin{proof}
Let $X\in \X(M)$ such that $AX=\pm X$. Then
\begin{equation*}
    A(JX)=-JAX=\mp JX \ .
\end{equation*}
Hence, $JX$ is in the $(\mp 1)$-eigen-subbundle of $A$. The same argument applies to the eigen-subbundles of $B$ 
and to the second statement.
\end{proof}

The $(\pm1)$-eigenbundles of $A$ play a special role in the Born structure, which will become apparent in 
Proposition~\ref{prop: almost Kunneth}. From now on, we will denote them by $L_{\pm}$. By 
Lemma~\ref{lem: J_interchanges_esp}, the endomorphism field $J$ defines isomorphisms
\begin{align*}
    J\colon L_{\pm}\longrightarrow L_{\mp} \ .
\end{align*}
In particular, $L_+$ and $L_-$ both have rank equal to $n=\frac{1}{2}\dim (M)$.

We will now study the transformation properties of $\omega$, $h$ and $g$ under the three recursion  operators. 

\begin{lemma}
\label{lem: omega_relations}
The 2-form $\omega$ transforms under the fields of endomorphisms $A$, $B$ and $J$ as
\begin{equation*}
\omega(JX,JY)=\omega(X,Y)=-\omega(AX,AY)=-\omega(BX,BY) \ .
\end{equation*}

\end{lemma}
\begin{proof}
For the first equality, we use the symmetry of $h$ and the anti-symmetry of $\omega$ to obtain
\begin{equation*}
    \omega(JX,JY)=-h(X, JY)= -h(JY,X)=-\omega(Y,X)=\omega(X,Y) \ .
\end{equation*}
Similarly, using the anti-symmetry of $\omega$ and the symmetry of $g$ we find
\begin{equation*}
    \omega(AX,AY)=g(X,AY)=g(AY,X)=\omega(A^2Y,X)=\omega(Y,X)=-\omega(X,Y) \ .
\end{equation*}
Finally, $B=-AJ$, and the first two equalities yield
\begin{equation*}
    \omega(BX,BY)=\omega(AJX,AJY)=-\omega(JX,JY)=-\omega(X,Y) \ .
\end{equation*}
This completes the proof.
\end{proof}

\begin{lemma}
 \label{lem: g transforms}
The pseudo-Riemannian metric $g$ satisfies
\begin{equation*}
    g(BX,BY)=g(X,Y)=-g(AX,AY)=-g(JX,JY) \ .
\end{equation*}
\end{lemma}
\begin{proof}
We apply Lemma~\ref{lem: omega_relations} several times and employ the fact that the recursion operators anti-commute. 
For the first equality, we have
\begin{equation*}
    g(BX,BY)=\omega(ABX,BY)=-\omega(BAX,BY)=\omega(AX,Y)=g(X,Y) \ .
\end{equation*}
Similarly, we proceed for $A$, obtaining
\begin{equation*}
    g(AX,AY)=\omega(X,AY)=-\omega(AX,Y)=-g(X,Y) \ .
\end{equation*}
For $J$, it follows that
\begin{equation*}
    g(JX,JY)=\omega(AJX,JY)=-\omega(JAX,JY)=-\omega(AX,Y)=-g(X,Y) \ .
\end{equation*}
\end{proof}

Proceeding analogously for $h$, we find that it is invariant under all three recursion operators:
\begin{lemma}
\label{lem: h transforms}
The pseudo-Riemannian metric $h$ is invariant under the three recursion operators, i.e.
\begin{equation*}
    h(AX,AY)=h(X,Y)=h(BX,BY)= h(JX,JY) \ .
\end{equation*}
\end{lemma}
In particular this shows that the pair $(h,J)$ forms an almost pseudo-Hermitian structure with fundamental
$2$-form $\omega$. This explains the choice of sign for $J$ in Definition~\ref{def: Born structure}. 

The transformation properties of the three bilinear forms are summarized in Table \ref{tab:comp_bil_forms}. 

\begin{table}[h]
    \centering
\begin{tabular}{ c|c|c } 
 $g(AX,AY)=-g(X,Y)$ & $h(AX,AY)=h(X,Y)$ & $\omega(AX,AY)=-\omega(X,Y)$\\
  $g(AX,Y)=-g(X,AY)$ & $h(AX,Y)=h(X,AY)$ & $\omega(AX,Y)=-\omega(X,AY)$\\
  &  &  \\
  $g(BX,BY)=g(X,Y)$ & $h(BX,BY)=h(X,Y)$ & $\omega(BX,BY)=-\omega(X,Y)$\\
  $g(BX,Y)=g(X,BY)$ & $h(BX,Y)=h(X,BY)$ & $\omega(BX,Y)=-\omega(X,BY)$\\
   & & \\
$g(JX,JY)=-g(X,Y)$ & $h(JX,JY)=h(X,Y)$ & $\omega(JX,JY)=\omega(X,Y)$ \\ 
 $g(JX,Y)=g(X,JY)$ & $h(JX,Y)=-h(X,JY)$ & $\omega(JX,Y)=-\omega(X,JY)$ \\ 
 \\
\end{tabular}
    \caption{Transformation properties of $g$, $h$ and $\omega$ under $A,B$ and $J$.}
    \label{tab:comp_bil_forms}
\end{table}

We will now prove further properties of the eigenbundles of $A$ and $B$ and the bilinear forms.

\begin{proposition}
\label{prop: almost Kunneth}
The $(\pm 1)$-eigenbundles of $A$ are Lagrangian for $\omega$. 
In particular, $(\omega, L_+,L_-)$ is an almost K\"unneth structure. 
The pseudo-Riemannian metric $g$ of the Born structure is the neutral metric 
associated with this almost K\"unneth structure.
\end{proposition}
\begin{proof}
Let $X,Y\in \X(M)$ such that $X$ and $Y$ are both in $L_+$ or both in $L_-$. Then by Lemma \ref{lem: omega_relations}, we have
\begin{equation*}
    \omega(X,Y)=\omega(AX,AY)=-\omega(X,Y).
\end{equation*}
Hence, the $(\pm 1)$-eigenbundles of $A$ are isotropic for $\omega$. Since they are complementary, they are also Lagrangian.

Note that the almost product structure $A$ agrees (up to sign) with the almost product structure $I$ associated to the almost K\"unneth structure that we defined in Subsection ~\ref{sect: def almost kunneth}. Since the neutral metric associated with the almost K\"unneth structure was defined by plugging $I$ into $\omega$, it follows that it is the same as the pseudo-Riemannian metric $g$ of the Born structure.
\end{proof}

The relation between almost K\"unneth structures and Born structures will be explored further in Section~\ref{sect: kunneth and born}.

\begin{corollary}
The signature of $g$ is neutral and the subbundles $L_{\pm}$ are null for $g$.
\end{corollary}
\begin{proof}
This is a reformulation of Proposition~\ref{prop: neutral metric}. The fact that the $L_{\pm}$ are null for $g$ can either be deduced from the transformation property of $g$ under $A$ or by observing that $g|_{L_{\pm}}=\pm\omega|_{L_{\pm}}$. Note that this holds for any neutral metric associated with an almost K\"unneth structure.
\end{proof}

\begin{proposition}\label{prop:orth}
The $(\pm 1)$-eigenbundles of $B$ are $g$-orthogonal.
\end{proposition}
\begin{proof}
Let $X,Y\in \X(M)$ such that $BX=X$, $BY=-Y$. Then
\begin{equation*}
    g(X,Y)=g(BX,BY)=-g(X,Y),
\end{equation*}
and hence $g(X,Y)=0$.
\end{proof}
Similarly, using Lemma \ref{lem: h transforms}, we find
\begin{proposition}
\label{prop: distr h-orthog}
The $(\pm 1)$-eigenbundles of $A$ and $B$ are $h$-orthogonal.
\end{proposition}

With the help of Proposition \ref{prop: distr h-orthog}, we can also say something about the signature of $h$:

\begin{proposition}
\label{prop: sign h}
The signature of $h$ is of the form $(2p,2q)$, where $p+q=n$.
\end{proposition}
\begin{proof}
As the subbundles $L_{\pm}$ are $h$-orthogonal, the signature of $h$ is the sum of the signatures of the restrictions $h|_{L_+}$ and $h|_{L_-}$. 
Since $J\colon L_+\to L_-$ is an isometry for $h$, both restricted metrics have the same signature $(p,q)$ and the assertion follows.
\end{proof}

We have now established the most basic properties of Born structures and are ready to construct first examples.
\begin{example}
\label{ex: cn}
Consider $\mathbb{C}^n$ with the standard Hermitian metric $h$ and the standard K\"ahler form $\omega$. The recursion operator $J$ from $h$ to $\omega$ is the standard complex structure on $\mathbb{C}^n$. Furthermore, thinking of $\mathbb{C}^n$ as $\mathbb{R}^{2n}$, let $L_{\pm}$ be the Lagrangian subbundles defined by the foliations given by the copies $\mathbb{R}^n\times\{y\}$ and $\{x\}\times\mathbb{R}^n$ of the two factors in the product decomposition $\mathbb{R}^{2n}=\mathbb{R}^n\times\mathbb{R}^n$. Then we define $A$ to be the almost product structure determined by $L_{\pm}$. We denote by $g$ the neutral metric defined by $\omega$ and $A$. Then $(g,h,\omega)$ is a Born structure on $\mathbb{C}^n$.
\end{example}

\begin{example}
\label{ex: tori pd}
Since the Born structure described in Example \ref{ex: cn} is invariant under translations, it descends to tori $\mathbb{C}^n/\Lambda$, where $\Lambda\subset\mathbb{C}^n$ is a lattice.
\end{example}

In view of Proposition \ref{prop: sign h} we can slightly modify Example \ref{ex: tori pd} to show that on tori $T^{2n}$, the signature of $h$ can take any form $(2p,2q)$ with $p+q=n$: 

\begin{example}
\label{ex: tori id}
Think of $T^{2n}$ as the product $T^{2p}\times T^{2q}$ and define $\omega$ and $g$ as in Example \ref{ex: tori pd}. However, for the almost complex structure, we choose now $J=J_p\oplus -J_q$, where $J_p$ and $J_q$ are the complex structures on the factors $T^{2p}$ and $T^{2q}$ coming from the standard complex structure on $\mathbb{C}^p$ and $\mathbb{C}^q$, respectively. The pseudo-Hermitian metric defined by $\omega$ and $J$ is now of signature $(2p,2q)$.
\end{example}

\subsection{Integrability of Born structures}
\label{sect: integrability}

 For integrability of a Born structure, there are several conditions that can be imposed. First of all, there is the condition of closedness of the $2$-form. Second of all, we could ask for integrability of the recursion operators. Since these all square to $\pm\id$, a necessary and sufficient condition for their integrability is the vanishing of the \emph{Nijenhuis tensor}.
 
 \begin{definition}
  The \emph{Nijenhuis tensor} of an endomorphism field $T$ is the $(1,2)$-tensor field defined by
 \begin{align*}
    N_T(X,Y)=[TX,TY]+T^2[X,Y]-T[TX,Y]-T[X,TY] \ .
 \end{align*}
 \end{definition}

That the vanishing of the Nijenhuis tensor is equivalent to the integrability of an almost complex structure is just the Newlander--Nirenberg theorem.  For an almost product structure $T$, the statement follows by observing that for $X,Y\in TM$, we have
\begin{align*}
    N_T(X,Y)=2\cdot\lb [X_+,Y_+]_-+[X_-,Y_-]_+\rb,
\end{align*}
where the subscripts $\pm$ denote the projections to the $(\pm 1)$-eigenbundles of $T$. 
Therefore, the vanishing of the Nijenhuis tensor is equivalent to the Frobenius integrability of the eigen-subbundles of $T$ and  
hence, to the integrability of $T$. 

The recursion operators of a Born structure have the special feature that the integrability of any two of them 
implies the integrability of the third one. 
 \begin{proposition}
 \label{prop: rec op}
 If two out of the three recursion operators in a Born structure are integrable, then so is the third one.
 \end{proposition}
\begin{proof}
For the Nijenhuis tensors of $A$, $B$ and $J$ a straightforward calculation leads to
 \begin{align*}
     N_J(X,Y)+N_J(AX,AY)=N_A(BX,BY)- N_A(X,Y) -A\lb N_B(AX,Y)+N_B(X,AY)\rb \ .
 \end{align*}
 If we assume $N_A\equiv N_B\equiv 0$, then this formula proves the vanishing of $N_J(X,Y)$ if $X$ and $Y$ 
 are in the same eigenbundle of $A$. For $X\in L_+$, $Y\in L_-$, we find
 \begin{align*}
     N_J(X,Y)=&[JX,JY]-[X,Y]-J[JX,Y]-J[X,JY]\\
    =& -[BX,BY]-[X,Y]-BA[BX,Y]+BA[X,BY]\\
    =&-N_B(X,Y)=0 \ .
 \end{align*}
Here, we used for the third equality that since $B$ interchanges the subbundles $L_{\pm}$ and $A$ is integrable, we have $[BX,Y]\in L_-$ and $[X,BY]\in L_+$. In the same way, one shows that $N_J\equiv N_A\equiv 0$ implies that $N_B$ vanishes identically. 

If we assume $N_J\equiv N_B\equiv 0$, then we obtain $N_A(X,Y)=N_A(BX,BY)$. It follows that $N_A(X,Y)$
vanishes if $X$ and $Y$ are in different eigenbundles for $B$. If $X$ and $Y$ are both in the $(+1)$-eigenbundle
of $B$, we compute
\begin{align*}
N_A(X,Y) &= [AX,AY]+[X,Y]- A[AX,Y]-A[X,AY]\\
&= -B[JX,JY]+B[X,Y]+BJ[JX,Y]+BJ[X,JY]\\
&=-BN_J(X,Y)=0 \ ,
\end{align*}
where we used $A=-JB=BJ$, and that the eigenbundles of $B$ are closed under commutators.
In the same way one proves that if $X$ and $Y$ are both in the $(-1)$-eigenbundle of $B$,
then
$$
N_A(X,Y)=BN_J(X,Y)=0 \ .
$$
This completes the proof.
 \end{proof}
 A statement equivalent to Proposition~\ref{prop: rec op} appeared in~\cite{zamkovoy2006geometry} with a different proof.
 
 \begin{remark}\label{rem:complexprod}
In the situation of Proposition~\ref{prop: rec op} the pair $(J,A)$ forms a complex product structure in the sense 
of Andrada and Salamon~\cite{AS}, cf.~Lemma~\ref{lem: para-quaternionic}. This is sometimes useful in order to exclude existence of integrable 
Born structures in the sense of the following definition, see Example~\ref{nonex} below.
\end{remark}
 
 \begin{definition}
 A Born structure is \emph{integrable}, if the two-form $\omega$ is closed and at least two, and therefore all, recursion 
 operators are integrable. 
  \end{definition}

We have already encountered first examples of integrable Born structures. 
The Born structures constructed on $\mathbb{C}^n$ in Example~\ref{ex: cn} and on tori in 
Examples~\ref{ex: tori pd} and \ref{ex: tori id} are clearly integrable.

Integrability of a Born structure implies that $(h,J)$ is a pseudo-K\"ahler structure with K\"ahler form $\omega$ and $(\omega,L_+,L_-)$ is a K\"unneth structure. 

\subsection{Born structures as enhanced almost K\"unneth structures}
\label{sect: kunneth and born}

We have seen in Proposition~\ref{prop: almost Kunneth} that every Born structure induces an almost K\"unneth structure. 
It turns out that the converse is also true and every almost K\"unneth structure can be realized as part of a Born structure.
Indeed, given an almost K\"unneth structure $(\omega,F,G)$, we can choose an isomorphism 
\begin{equation*}
    \Tilde{J}\colon F\to G \ , 
\end{equation*}
satisfying
\begin{equation*}
    \omega(\Tilde{J}X,Y)=-\omega(X,\Tilde{J}Y),\quad X,Y\in F \ .
\end{equation*}
Locally, this can be realized as follows. Choose a local frame $\{f_i,g_j\}_{i,j=1}^n$ such that the $\{f_i\}$ and $\{g_i\}$ 
are a local frame for $F$ and $G$, respectively, and $\omega(f_i,g_j)=\delta_{ij}$. Then set $\tilde{J}(f_i)=g_i$.

Using $\tilde{J}$, we can define an endomorphism field $J$ on $TM$ by setting
\begin{equation*}
    J\mid_{F}:=\Tilde{J},\quad J\mid_{G}:=-\Tilde{J}^{-1}.
\end{equation*}
By construction, $J$ is an almost complex structure which interchanges the Lagrangian subbundles of the almost 
K\"unneth structure and such that $h(X,Y)=\omega(X,JY)$ defines a pseudo-Riemannian metric $h$. 
Then we obtain a Born structure $(g,h,\omega)$, where $g$ is the neutral metric associated with the almost K\"unneth structure. 
This makes sense even though $g$ is only well-defined up to sign. Choosing the neutral metric with a different sign still yields a Born structure $(-g,h,\omega)$, where the recursion operators $A$ and $B$ are replaced by $-A$ and $-B$. 

It follows that a Born structure is an almost K\"unneth structure, together with a choice of an almost complex structure $J$, 
which is compatible with the almost K\"unneth structure in the sense that it interchanges the Lagrangian subbundles and 
it satisfies $\omega(J\cdot,\cdot)=-\omega(\cdot,J\cdot)$. In particular, an integrable Born structure is a K\"unneth 
structure together with a compatible complex structure.

\begin{remark}
Note that if $J$ is defined by an isomorphism $\tilde{J}$, which is locally of the form described above, the resulting $h$ will always be positive definite. The reason is that the condition $\omega(f_i,g_j)=\delta_{ij}$ implies that $\{f_i\}$ is a local orthonormal frame of $L_+$ with respect to the $h$ so constructed. To obtain a Born structure with an indefinite $h$ in an analogous construction, we need a local frame $\{f_i,g_i\}_{i=1}^n$ such that $\omega(f_i,g_j)=\varepsilon_i\delta_{ij}$, where $\varepsilon_i=\pm 1$. However, in general this will not lead to a globally defined pseudo-Riemannian metric.
\end{remark}

\subsection{Relation to hypersymplectic geometry}\label{ss:hyper}

According to~\cite{bande2008geometry}, hypersymplectic structures in the sense of Hitchin~\cite{hitchin1990hypersymplectic}
can be defined by a diagram of symplectic forms 
\begin{center}
        \begin{tikzcd}[column sep=small]& \omega\arrow[dl,"A"']\arrow{dr}{B} & \\\alpha \arrow{rr}{J} & & \beta
     \end{tikzcd}
    \end{center}
such that the recursion operators satisfy $-J^2=A^2=B^2=\id$. This looks very similar to our definition of a Born structure, 
the only difference being that the symplectic forms $\alpha$ and $\beta$ are replaced by pseudo-Riemannian metrics. 

It was proved in~\cite[Section~8.3]{kuennethgeom} that every hypersymplectic structure 
induces an $S^1$-family of K\"unneth structures. The hypersymplectic metric from Hitchin's definition is in fact the neutral
metric $g$ associated with any one of these K\"unneth structures.
By Section~\ref{sect: kunneth and born} above, every K\"unneth 
structure in this family gives rise to a Born structure. In some situations, this even yields an $S^1$-family of integrable Born structures:
\begin{theorem}
\label{thm: s1-fam int bs}
Let $M$ be a smooth manifold admitting a hypersymplectic structure $(\omega,\alpha,\beta)$ with hypersymplectic metric $g$. 
Assume that there is an almost complex structure $\tilde{J}$ on $M$, which anti-commutes with  $A$ and $B$ and such that 
$g(\tilde{J}X,\tilde{J}Y)=-g(X,Y)$ for all $X$, $Y$ in $TM$. Then
\begin{center}
    \begin{tikzcd}[column sep=small]& g\arrow[dl,"I_{\theta}"']\arrow{dr}{\tilde{B}_{\theta}} & \\\beta_{\theta} \arrow{rr}{-\tilde{J}} & & h_{\theta}
     \end{tikzcd}
\end{center}
defines an $S^1$-family of Born structures, where
\begin{align*}
    \beta_{\theta}=-\sin(\theta)\alpha + \cos(\theta)\beta, \quad  I_{\theta}=\cos(\theta)A + \sin(\theta)B,
\end{align*}
and $\tilde{B}_{\theta}=\tilde{J}I_{\theta}$.
In particular, if $\tilde{J}$ is integrable, then every Born structure in this family is integrable.
\end{theorem}
\begin{proof}
Recall from~\cite[Section~8.3]{kuennethgeom} the $S^1$-family of K\"unneth structures $(\alpha_{\theta+\pi/2},\mathcal{F}_{\theta},\mathcal{G}_{\theta})$ that arises from the hypersymplectic structure,
where 
\begin{align*}
    \alpha_{\theta}=\cos(\theta)\alpha + \sin(\theta)\beta, 
\end{align*}
and $\mathcal{F}_{\theta},\mathcal{G}_{\theta}$ are the eigen-foliations of the product structure
\begin{align*}
  I_{\theta}=\cos(\theta)A + \sin(\theta)B.  
\end{align*}

Since $\tilde{J}$ anti-commutes with $A$ and $B$, it anti-commutes with $I_{\theta}$ for each $\theta$ and since every hypersymplectic structure in the family has the same associated neutral metric $g$ and $g(\tilde{J}X,Y)=g(X,\tilde{J}Y)$ by assumption, it follows that $\alpha_{\theta}(\tilde{J}X,Y)=-\alpha_{\theta}(X,\tilde{J}Y)$ for each $\theta$. Therefore, every K\"unneth structure in the $S^1$-family is compatible with $\tilde{J}$ in the sense of Section~\ref{sect: kunneth and born} and we obtain an $S^1$-family of Born structures
\begin{center}
        \begin{tikzcd}[column sep=small]& g\arrow[dl,"I_{\theta}"']\arrow{dr}{\tilde{B}_{\theta}} & \\\beta_{\theta} \arrow{rr}{-\tilde{J}} & & h_{\theta}
     \end{tikzcd},
\end{center}
where $\beta_{\theta}:=\alpha_{\theta+\pi/2}$ and $\tilde{B}_{\theta}=\tilde{J}I_{\theta}$.

If $\tilde{J}$ is integrable, then, since $\beta_{\theta}$ is closed and $I_{\theta}$ is integrable for each $\theta$, every Born structure in this family is completely integrable because of Proposition~\ref{prop: rec op}.
\end{proof}

We will discuss an example of this construction on $Nil^3\times \mathbb{R}$ in Section~\ref{s:examples} below.

\subsection{Comparison to the previous literature}
\label{sect: def born physics}

In~\cite{freidel2019unique} a Born structure was defined as an almost para-Hermitian structures $(g,I)$ with corresponding 
$2$-form $\omega$, together with a pseudo-Riemannian metric $h$ such that the recursion operator from $g$ to $h$ squares to $\id$ and the one from $\omega$ to $h$ squares to $-\id$. From Proposition~\ref{prop: almost Kunneth} it follows that this is equivalent to our Definition~\ref{def: Born structure}.

However, our notation differs from the notation that has been used in the physics literature. In \cite{freidel2019unique}, \cite{svoboda2020born} and \cite{marotta2021born}, the metric $h$ is denoted by $\mathcal{H}$ and the metric $g$ is called $\eta$. Moreover, the almost product structures $A$ and $B$ are denoted by $K$ and $J$, while the almost complex structure $J$ is called $I$. 

The structure formed by $A$, $B$ and $J$ is called an \emph{almost para-quaternionic structure} and was first introduced in \cite{libermann1952structures} and later studied in~\cite{vukmirovic2003quaternionic}, \cite{ivanov2005parahermitian}, \cite{zamkovoy2006geometry}, and in many other papers.

\section{Connections associated with a Born structure}
\label{sect: connections}

A Born structure on a smooth manifold $M$ determines the following affine connections:
\begin{itemize}
    \item $\cg$, the \emph{Levi-Civita connection} of the pseudo-Riemannian metric $g$,
    \item $\ch$, the \emph{Levi-Civita connection} of the pseudo-Riemannian metric $h$,
    \item $\ck$, the \emph{K\"unneth connection} of the almost K\"unneth structure $(\omega,L_+,L_-)$,
    \item $\cc$, the \emph{canonical connection} of the almost K\"unneth structure $(\omega,L_+,L_-)$,
    \item $\cb$, the \emph{Born connection}, which is compatible with $\omega,g$ and $h$ and has vanishing \emph{generalized torsion}.
\end{itemize}
We already discussed the K\"unneth connection in Section~\ref{sect: def almost kunneth}. We will now define the canonical connection and the Born connection and discuss their properties. Furthermore, we will explain the relations between the different connections.

\subsection{The canonical connection}\label{sect: cc}

Given an almost K\"unneth structure $(\omega,L_+,L_-)$ with almost product structure $A$, whose eigenbundles are 
$L_{\pm}$, we have the associated neutral Riemannian metric $g$ defined in Proposition~\ref{prop: neutral metric}.
Its Levi-Civita connection $\cg$ does not necessarily commute with $A$, and this failure motivates the definition of 
the canonical connection.

For $X,Y\in\X(M)$ the \emph{canonical connection} $\cc$ is defined by averaging $\cg$ under conjugation with $A$:
\begin{equation}\label{eq:cc}
 \cc_X Y= \frac{1}{2}\left(\cg_X Y + A\cg_X AY\right) \ .
\end{equation}
This does indeed define a connection which commutes with $A$. Furthermore, this new connection is still 
compatible with $g$, but may have non-trivial torsion. 
Since the canonical connection is compatible with $g$ and commutes with $A$, it is also compatible with $\omega$.
This means that the whole almost K\"unneth structure is parallel with respect to $\cc$.

One can rewrite the definition of the canonical connection as
\begin{equation*}
    \cc_X Y=\left(\cg_X Y_+\right)_+ +\left(\cg_X Y_-\right)_-,
\end{equation*}
where the subscripts $\pm$ denote the projections to the subbundles $L_{\pm}$. 

It follows immediately from this and the fact that the subbundles $L_{\pm}$ are null for $g$ that
\begin{align*}
g\lb\cc_X Y_+,Z_-\rb + g\lb Y_+, \cc_X Z_-\rb =g\lb\cg_X Y_+,Z_-\rb + g\lb Y_+, \cg_X Z_-\rb = L_X\lb g(Y_+,Z_-)\rb  .   
\end{align*}
Moreover, using again that the $L_{\pm}$ are null for $g$, we find
\begin{align*}
  g\lb\cc_X Y_{\pm},Z_{\pm}\rb + g\lb Y_{\pm},\cc_XZ_{\pm}\rb=0
\end{align*}
and 
\begin{align*}
L_X\lb g \lb Y_{\pm},Z_{\pm}\rb \rb=0 \ .
\end{align*}

We have seen in Section~\ref{sect: def almost kunneth}, that for the K\"unneth connection, just like the for canonical connection,  the 
whole almost K\"unneth structure is parallel. Therefore, the question arises whether there are conditions under which these connections coincide with each other,  and how they are related in general. By the uniqueness of the K\"unneth connection, it is clear that a necessary and sufficient condition for both connections to agree with each other is the vanishing of the mixed torsion of the canonical connection. For the purpose of stating the relation between the connections more precisely, we will need the following lemma:
\begin{lemma}
\label{lem: cg and omega}
The following equalities hold for $X,Y,Z\in\X(M)$:
\begin{enumerate}
    \item[i)] $\cg_X\omega\left(Y_{\pm},Z_{\mp}\right)=0$,
    \item[ii)] $\cg_X\omega\left(Y_{\pm},Z_{\pm}\right)=-2\omega\left(\cg_X Y_{\pm},Z_{\pm}\right)=-2\omega\left(Y_{\pm},\cg_X Z_{\pm}\right)$.
\end{enumerate}
\end{lemma}
\begin{proof}
The first statement follows from
\begin{align*}
 \cg_Z\omega(X_+,Y_-)&=L_Z\left(\omega(X_+,Y_-)\right)-\omega(\cg_Z X_+,Y_-)-\omega(X_+,\cg_Z Y_-)\\
    &=L_Z\left(\omega(X_+,Y_-)\right)-\omega(\cc_Z X_+,Y_-)-\omega(X_+,\cc_Z Y_-)\\
    &=0,
\end{align*}
where in the second line we used the definition of $\cc$ and in the last step that $\omega$ is parallel for $\cc$.

For the second statement, we observe
\begin{align*}
    \omega\left(Y_{\pm},\cg_X Z_{\pm}\right)=\pm g\left(Y_{\pm},\cg_X Z_{\pm}\right)=\mp g\left(\cg_X Y_{\pm}, Z_{\pm}\right)=\omega\left(\cg_X Y_{\pm}, Z_{\pm}\right),
\end{align*}
where we used for the second equality that the $L_{\pm}$ are null for $g$. This yields
\begin{align*}
  \cg_X\omega\left(Y_{\pm},Z_{\pm}\right)=-\omega\left(\cg_X Y_{\pm},Z_{\pm}\right)-\omega\left( Y_{\pm},\cg_X Z_{\pm}\right) = -2\omega\left(\cg_X Y_{\pm},Z_{\pm}\right).
\end{align*}
\end{proof}
With the help of Lemma \ref{lem: cg and omega}, we obtain the following result:
\begin{proposition}
\label{prop: omega_K}
The K\"unneth connection and the canonical connection are related by
\begin{align*}
    \omega\lb\ck_X Y,Z\rb=\omega\lb \cc_X Y,Z\rb - \Omega^K(X,Y,AZ),
\end{align*}
where
\begin{align*}
   \Omega^K(X,Y,AZ)=\frac{1}{2}\,\left\{\dd\omega(X,Y_+,Z_-)+\dd\omega(X,Y_-,Z_+)\right\}.
\end{align*}
\end{proposition}
\begin{proof}
Using the definition of the K\"unneth connection, we write
\begin{align*}
    \omega &\left(\ck_X Y_+,Z_-\right)= L_{X_+}\omega\left(Y_+,Z_-\right)-\omega\left(Y_+,[X_+,Z_-]\right)+\omega\left([X_-,Y_+],Z_-\right)\\
    &= \omega\lb\cc_{X_+}Y_+,Z_-\rb + \omega\lb Y_+, \cc_{X_+}Z_-\rb -\omega\left(Y_+,[X_+,Z_-]\right)+\omega\left([X_-,Y_+],Z_-\right),
\end{align*}
where we used for the second equality that the canonical connection is compatible with $\omega$. Using the definition of $\cc$ and that $\cg$ is torsion-free, we observe that
\begin{align*}
  \omega\lb Y_+, \cc_{X_+}Z_- -[X_+,Z_-]\rb= \omega\lb Y_+, \cg_{X_+}Z_- -[X_+,Z_-]\rb = \omega\lb Y_+, \cg_{Z_-}X_+\rb,
\end{align*}
and similarly,
\begin{align*}
    \omega\left([X_-,Y_+],Z_-\right)=\omega\lb \cg_{X_-}Y_+-\cg_{Y_+}X_-,Z_-\rb = \omega\lb\cc_{X_-}Y_+,Z_-\rb-\omega\lb\cg_{Y_+}X_-,Z_-\rb.
\end{align*}
It follows that
\begin{align*}
    \omega\left(\ck_X Y_+,Z\right)=\omega\lb \cc_X Y_+,Z\rb +  \omega\lb Y_+, \cg_{Z_-}X_+\rb-\omega\lb\cg_{Y_+}X_-,Z_-\rb.
\end{align*}
By analogous arguments we obtain
\begin{align*}
    \omega\left(\ck_X Y_-,Z\right)=\omega\lb \cc_X Y_-,Z\rb +  \omega\lb Y_-, \cg_{Z_+}X_-\rb-\omega\lb\cg_{Y_-}X_+,Z_+\rb.
\end{align*}
Therefore, we have
\begin{align*}
    \omega\lb\ck_X Y,Z\rb=\omega\lb \cc_X Y,Z\rb - \Omega^K(X,Y,AZ),
\end{align*}
where
\begin{align*}
   \Omega^K(X,Y,AZ)=&-\omega\lb Y_+, \cg_{Z_-}X_+\rb+\omega\lb\cg_{Y_+}X_-,Z_-\rb\\ &- \omega\lb Y_-, \cg_{Z_+}X_-\rb+\omega\lb\cg_{Y_-}X_+,Z_+\rb.
\end{align*}
Using Lemma \ref{lem: cg and omega} twice, we find
\begin{align*}
\omega\lb Y_+, \cg_{Z_-}X_+\rb &= \frac{1}{2}\,\cg_{Z_-}\omega\lb X_+,Y_+\rb\\
&= \frac{1}{2}\,\left\{ \cg_{X_+}\omega\lb Y_+,Z_-\rb -\cg_{Y_+}\omega\lb Z_-,X_+\rb + \cg_{Z_-}\omega\lb X_+,Y_+\rb \right\}\\
&= \frac{1}{2}\,\dd\omega(X_+,Y_+,Z_-).
\end{align*}
Proceeding similarly with the other terms, it follows that
\begin{align*}
   \Omega^K(X,Y,AZ)=&\frac{1}{2}\,\dd\omega(X_+,Y_+,Z_-)+\frac{1}{2}\,\dd\omega(X_-,Y_+,Z_-)\\ &+\frac{1}{2}\,\dd\omega(X_-,Y_-,Z_+)+\frac{1}{2}\,\dd\omega(X_+,Y_-,Z_+)\\
   =& \frac{1}{2}\,\dd\omega(X,Y_+,Z_-)+\frac{1}{2}\,\dd\omega(X,Y_-,Z_+).
\end{align*}
\end{proof}

The subbundles $L_{\pm}$ induce a bigrading on the differential forms of $M$. In particular, if the subbundles are integrable, then the exterior derivative splits as $d=d_+ + d_-$, where $d_+$ and $d_-$ are of bidegree $(1,0)$ and $(0,1)$, respectively. By Proposition \ref{prop: almost Kunneth}, the 2-form $\omega$ is of type $(1,1)$. Using this bigrading, we obtain the following necessary and sufficient conditions for the vanishing of the difference of the K\"unneth connection and the canonical connection:

\begin{corollary}
The K\"unneth connection and the canonical connection agree with each other if and only if
\begin{align*}
    \dd\omega^{(1,2)}=0=\dd\omega^{(2,1)} \ .
\end{align*}
In particular, if $\omega$ is closed, then $\cc=\ck$. If $A$ is integrable, then $\cc=\ck$ if and only if $\omega$ is closed.
\end{corollary}
Using Theorem \ref{thm: ck is LC}, this implies:
\begin{corollary}
\label{cor: Kunneth cc}
If $(\omega,L_+,L_-)$ is K\"unneth, then $\cc=\ck=\cg$.
\end{corollary}

\subsection{The Born connection}\label{sect:def_cb}

We will now discuss connections compatible with the full Born structure. 

\begin{definition}
A connection $\nabla$ is \emph{compatible with the Born structure} $(g,h,\omega)$ if
\begin{equation*}
    \nabla h=\nabla\omega=\nabla g=0.
\end{equation*}
\end{definition}

First we show that compatible connections exist, using again the averaging mechanism used in the 
definition of the canonical connection.
\begin{proposition}\label{prop:Bornconnection}
Let $(g,h,\omega)$ be a Born structure, and $\ck$ the K\"unneth connection of the 
underlying almost K\"unneth structure. Then 
\begin{equation*}
    \nabla_X Y= \frac{1}{2}\left(\ck_X Y + B\ck_X BY\right)  = \frac{1}{2}\left(\ck_X Y- J\ck_X JY\right) 
\end{equation*} 
defines a connection that is compatible with the Born structure.
\end{proposition}
\begin{proof}
We first check that the two expressions involving $B$ and $J$ respectively agree. 
For this we use that $B=JA$, that $J$ and $A$ anti-commute, 
and the fact that $\ck$
commutes with $A$. Together, these identities imply
$$
B\ck_X BY = JA\ck_X JAY = JA\ck_X (-AJY)= -JA^2\ck_X JY = -J\ck_X JY \ .
$$
Thus the two expressions do indeed agree. They define an affine connection $\nabla$,
for which we want to show that it is compatible with the  Born structure. 

Recall that the K\"unneth connection $\ck$ is compatible with $g$ and $\omega$. It commutes with $A$, but not necessarily with $B$.
However, $\nabla$ does commute with $B$:
\begin{alignat*}{1}
 \nabla_X BY &= \frac{1}{2}\left(\ck_X BY + B\ck_X B^2Y\right)\\
 &= \frac{1}{2}\left(\ck_X BY + B\ck_X Y\right)\\
 &= B\left(\frac{1}{2}\left(B\ck_X BY + \ck_X Y\right)\right)\\
 &=B\nabla_X Y \ .
\end{alignat*}
Moreover, $\nabla$ is still compatible with $g$ and $\omega$, and commutes with $A$. Since it also commutes with $B$ it
commutes with $J$ as well, and is also compatible with $h$.
\end{proof}

In general there are many connections compatible with a Born structure, but 
from our point of view, the connection constructed above is the most natural,
and would deserve to be called the Born connection. However, this name is already in use for a connection that
does not always agree with this $\nabla$.

Freidel, Rudolph and Svoboda~\cite{freidel2019unique}  used a condition analogous to torsion-freeness
in the case of the Levi-Civita connection of a pseudo-Riemannian metric in order to single out a preferred 
connection in the space of all connections compatible with a Born structure. Their condition is the 
following:
 \begin{definition}
 Let $\cc$ be the canonical connection for an almost K\"unneth structure with associated neutral metric $g$ on a smooth manifold $M$. Then an affine connection $\nabla$ on $TM$ has \emph{vanishing generalized torsion} if
 \begin{equation}
 \label{eq: gen torsion}
    g(\nabla_X Y-\nabla_Y X,Z)+g(\nabla_Z X,Y)=g(\cc_X Y-\cc_Y X,Z)+g(\cc_Z X,Y)
\end{equation}
for all $X,Y,Z\in\X(M)$.
 \end{definition}
This definition does not arise naturally in our setup, and we refer to~\cite{freidel2019unique} for its motivation.

The following was proved in \cite{freidel2019unique}:
\begin{theorem}\cite[Theorem~1]{freidel2019unique}\label{unique}
Given a Born structure $(g,h,\omega)$, there is a unique connection, the \emph{Born connection} $\cb$, that is compatible with the Born structure and has vanishing generalized torsion.
\end{theorem}

In the case of an integrable Born structure the Born connection does agree with our candidate $\nabla$
constructed by averaging the K\"unneth connection under the conjugation with $B$ or $J$, although this is 
far from obvious.
\begin{theorem}\label{main}
Let $(g,h,\omega)$ be an integrable Born structure. Then the connection $\nabla$ defined 
in Proposition~\ref{prop:Bornconnection} agrees with the Born connection $\cb$.
\end{theorem}

\begin{proof}
By Proposition~\ref{prop:Bornconnection}, $\nabla$ is compatible with the Born structure. Therefore,
by Theorem~\ref{unique}, it suffices to show that $\nabla$ has vanishing generalized torsion, and for this we
use integrability.

From Corollary~\ref{cor: Kunneth cc} we know that in the integrable case the canonical and K\"unneth 
connections coincide. Therefore, on the right-hand side of~\eqref{eq: gen torsion} we can use the K\"unneth
connection. Thus, using the definition of $\nabla$ in terms of $\ck$, the condition~\eqref{eq: gen torsion} is equivalent to 
 \begin{equation}
 \label{eq: gen torsion2}
    g(\ck_X Y-\ck_Y X,Z)+g(\ck_Z X,Y)=g(B\ck_X BY-B\ck_Y BX,Z)+g(B\ck_Z BX,Y) \ .
\end{equation}
Since in the integrable case $\ck$ is torsion-free in the usual sense, the first summand on the left-hand side 
of~\eqref{eq: gen torsion2} simplifies to $g([X,Y],Z)$. Now the integrability of $B$ implies 
$$
[X,Y] = B[BX,Y]+B[X,BY]-[BX,BY] \ ,
$$
where on the right-hand side we can rewrite commutators in terms of $\ck$ since this is torsion-free:
$$
[X,Y] = B\ck_{BX}Y-B\ck_{Y}BX+B\ck_X BY - B\ck_{BY} X -[BX,BY] \ .
$$
Substituting this into the left-hand side of~\eqref{eq: gen torsion2} we see that certain summands on
the two sides agree and therefore drop out. We are left with the condition
\begin{equation}
 \label{eq: gen torsion3}
    g(B\ck_{BX}Y-B\ck_{BY} X-[BX,BY],Z)+g(\ck_Z X,Y)=g(B\ck_Z BX,Y) \ .
\end{equation}
\begin{lemma}
The identity~\eqref{eq: gen torsion3} holds if $X$ and $Y$ are in the same eigenbundle for $B$.
\end{lemma}
\begin{proof}
Suppose $BX=X$ and $BY=Y$. Then~\eqref{eq: gen torsion3} becomes
\begin{equation}
 \label{eq: gen torsion4}
    g(B\ck_{X}Y-B\ck_{Y} X-[X,Y],Z)+g(\ck_Z X,Y)=g(B\ck_Z X,Y) \ .
\end{equation}
Since the eigenbundles of $B$ are integrable, we have $[X,Y]=B[X,Y]$, and so the first summand on 
the left-hand side vanishes by the torsion-freeness of $\ck$. For the second summand we have:
$$
g(\ck_Z X,Y) = g(B\ck_Z X,BY)  = g(B\ck_Z X,Y)   \ .
$$
This show that~\eqref{eq: gen torsion3} holds if $BX=X$ and $BY=Y$. The proof also works if $BX=-X$ and $BY=-Y$.
\end{proof}

\begin{lemma}
The identity~\eqref{eq: gen torsion3} holds if $X$ and $Y$ are in different eigenbundles for $B$.
\end{lemma}
\begin{proof}
Suppose $BX=X$ and $BY=-Y$. Then~\eqref{eq: gen torsion3} becomes
\begin{equation}
 \label{eq: gen torsion5}
    g(B\ck_{X}Y+B\ck_{Y} X+[X,Y],Z)+g(\ck_Z X,Y)=g(B\ck_Z X,Y) \ .
\end{equation}
The right-hand side is
$$
g(B\ck_Z X,Y) = g(\ck_Z X,BY) = -g(\ck_Z X,Y) \ ,
$$
which is the negative of the second summand on the left-hand side. The first summand there 
can be rewritten using once more the torsion-freeness of $\ck$:
\begin{equation}
 \label{eq: gen torsion6}
\begin{split}
g(B\ck_{X}Y+B\ck_{Y} X+[X,Y],Z) &= g(B\ck_{X}Y+B\ck_{Y} X+\ck_X Y - \ck_Y X,Z)\\
&=g((B\ck_{X}Y+\ck_X Y ) +(B\ck_{Y} X- \ck_Y X),Z) \ .
\end{split}
\end{equation}
Note that $B\ck_{X}Y+\ck_X Y$ is in the $(+1)$-eigenbundle of $B$, and $B\ck_{Y} X- \ck_Y X$ 
is in the $(-1)$-eigenbundle. These eigenbundles are $g$-orthogonal by Proposition~\ref{prop:orth}, and therefore
we now consider separately the two cases where $Z$ is in one of them.

If $BZ = Z$, then using~\eqref{eq: gen torsion6} the first summand on the left-hand side of~\eqref{eq: gen torsion5}
is
\begin{equation*}
\begin{split}
2g(\ck_X Y,Z) &= 2L_X g(Y,Z)-2g(Y,\ck_X Z)\\
&= -2g(Y,\ck_Z X+[X,Z])\\
&=-2g(Y,\ck_Z X) \ ,
\end{split}
\end{equation*}
which shows that~\eqref{eq: gen torsion5} holds in this case. In this calculation we used that $\ck$ is 
compatible with $g$, that $[X,Z]$ is in the $(+1)$-eigenbundle of $B$ by the integrability of this subbundle,
and that this is $g$-orthogonal to the $(-1)$-eigenbundle containing $Y$.

If $BZ = -Z$, then using~\eqref{eq: gen torsion6} the first summand on the left-hand side of~\eqref{eq: gen torsion5}
is
\begin{equation*}
\begin{split}
-2g(\ck_Y X,Z) &= -2L_Y g(X,Z)+2g(X,\ck_Y Z)\\
&= 2g(X,\ck_Z Y+[Y,Z])\\
&=2g(X,\ck_Z Y) \\
&= 2L_Z g(X,Y)-2g(\ck_Z X,Y)\\
&=-2g(\ck_Z X,Y) \ ,
\end{split}
\end{equation*}
where we use the same arguments as before. This shows that~\eqref{eq: gen torsion5} holds in the case $BZ = -Z$.
By linearity it holds for all $Z$.

In exactly the same way one proves the case when $BX=-X$ and $BY=Y$.
\end{proof}
Combining the two lemmas with linearity, we see that the condition~\eqref{eq: gen torsion3} holds for all $X$, $Y$ and
$Z$. This finally completes the proof of the theorem.
\end{proof}

\begin{example}
Consider the situation of Theorem~\ref{thm: s1-fam int bs}, where an $S^1$-family of Born structures is 
constructed from a hypersymplectic structure. In that case, the underlying K\"unneth structures all have 
for their K\"unneth connection $\ck$ the Levi-Civita connection $\cg$ of the hypersymplectic metric. Moreover,
the almost complex structure $\tilde{J}$ is independent of the parameter $\theta\in S^1$. Therefore, the 
connection $\nabla$ constructed in Proposition~\ref{prop:Bornconnection} by averaging $\ck$ under the 
conjugation with the almost complex structure is also independent of $\theta$. Whenever $\tilde{J}$ is 
integrable, this is the Born connection $\cb$, which is then the same for all Born structures in the $S^1$-family.
\end{example}

Finally we note that in the special case where the K\"unneth connection $\ck$ commutes with $B$, the 
connection $\nabla$ we defined equals the K\"unneth connection. Therefore, Theorem~\ref{main} has the 
following consequence.
\begin{corollary}
For an integrable Born structure with $B$ or, equivalently, $J$, parallel with respect to $\ck$, we have $\cc = \cg=\ck=\cb$.
\end{corollary}
In this special situation $\cb$ is torsion-free in the usual sense. 
In general, the torsion of $\cb$ measures the failure of $\ck$ to commute with $B$. This is the 
content of the following proposition.
\begin{proposition}
Let $(g,h,\omega)$ be an integrable Born structure, and $T^B$ the torsion tensor of the Born connection
$\cb$. 
Then $T^B(X,Y)$ vanishes if $X$ and $Y$ are in the same eigenbundle of $B$. 

If $BX=X$ and $BY=-Y$, then 
\begin{alignat*}{1}
T^B(X,Y) &=- \frac{1}{2}\left(\ck_X Y+B\ck_X Y\right)+\frac{1}{2}\left(\ck_Y X - B\ck_Y X\right)\\
&=-\pi_+ ( \ck_X Y)+\pi_- ( \ck_Y X)  \ ,
\end{alignat*}
where $\pi_{\pm}$ are the projections to the $(\pm 1)$-eigenbundles of $B$.

In particular, $T^B$ vanishes identically if and only if $\ck$ commutes with $B$.
\end{proposition}
\begin{proof}
Using Theorem~\ref{main}, we shall calculate with the formula from Proposition~\ref{prop:Bornconnection}.

Assume first that $BX=X$ and $BY=Y$. Then 
\begin{alignat*}{1}
T^B(X,Y) &= \cb_X Y - \cb_Y X - [X,Y] \\
&=\frac{1}{2}\left(\ck_X Y + B\ck_X BY\right) -\frac{1}{2}\left(\ck_Y X + B\ck_Y BX\right) -[X,Y]\\
&=\frac{1}{2}\left(\ck_X Y + B\ck_X Y\right) -\frac{1}{2}\left(\ck_Y X + B\ck_Y X\right) -[X,Y]\\
&=\pi_+ (\ck_X Y )-\pi_+ (\ck_Y X) - \pi_+ ([X,Y])\\
&=\pi_+ (\ck_X Y-\ck_Y X - [X,Y])\\
&=\pi_+ (0) = 0 \ ,
\end{alignat*}
where we have used $B[X,Y]=[X,Y]$ by the assumptions on $X$ and $Y$ and the integrability of $B$,
and the torsion-freeness of $\ck$ in the integrable case.
The same argument works if $BX=-X$ and $BY=-Y$.

Next assume that $BX=X$ and $BY=-Y$. The torsion-freeness of $\ck$ implies
\begin{alignat*}{1}
T^B(X,Y) &=\frac{1}{2}\left(\ck_X Y + B\ck_X BY\right) -\frac{1}{2}\left(\ck_Y X + B\ck_Y BX\right) -[X,Y]\\
&=\frac{1}{2}\left(\ck_X Y - B\ck_X Y\right) -\frac{1}{2}\left(\ck_Y X + B\ck_Y X\right) -\ck_X Y +\ck_Y X\\
&=-\frac{1}{2}\left(\ck_X Y + B\ck_X Y\right) +\frac{1}{2}\left(\ck_Y X - B\ck_Y X\right)\\
&=-\pi_+ ( \ck_X Y)+\pi_- ( \ck_Y X)  \ ,
\end{alignat*}
as claimed.
\end{proof}

\section{Examples of integrable Born structures}\label{s:examples}

In this section we want to provide some examples of integrable Born structures on closed manifolds that go beyond the 
rather obvious ones we have seen in Examples~\ref{ex: tori pd} and \ref{ex: tori id}.

Given an integrable Born structure
    \begin{center}
        \begin{tikzcd}[column sep=small]& g\arrow[dl,"A"']\arrow{dr}{B} & \\\omega \arrow{rr}{-J}   & &h \ \ \ ,
     \end{tikzcd}
    \end{center}
we know that $(h,J)$ is a pseudo-K\"ahler structure with K\"ahler form $\omega$, and $(\omega,L_+,L_-)$ is a K\"unneth structure. 
Therefore, we will look for examples in classes of manifolds for which it is known that one or both of these structures occur.

A very tractable class of manifolds consists of nilmanifolds. Note that for a simply connected nilpotent Lie group $G$ admitting a lattice 
$\Gamma$ any left-invariant integrable Born structure descends to the compact nilmanifold $\Gamma\setminus G$. Therefore, to investigate 
whether a nilmanifold is left-invariant Born, it suffices to work at the level of the Lie algebra. 

The definition of a Born structure can be transcribed to Lie algebras in the obvious way. We will call a Born structure on a Lie algebra 
$\mathfrak{g}$ \emph{integrable}, if its $2$-form is closed under the Chevalley--Eilenberg differential, the endomorphism 
$J$ has vanishing Nijenhuis tensor and the eigenspaces of $A$ are subalgebras of $\mathfrak{g}$. 
Furthermore, a Lie algebra admitting an integrable Born structure will be called \emph{Born}. 

Integrable Born structures on a Lie algebra yield left-invariant integrable Born structures on the corresponding Lie group. 
Since all of the nilpotent Lie algebras that we consider have a basis with rational structure constants, the corresponding 
Lie groups admit lattices and therefore, the integrable Born structures we find give rise to examples of compact Born manifolds.

A \emph{compatible pair} on a Lie algebra $\mathfrak{g}$ is a symplectic form $\Omega$ on $\mathfrak{g}$ together with a complex structure 
$J$ on $\mathfrak{g}$, such that $\Omega(JX,JY)=\Omega(X,Y)$ and, therefore, $h(\cdot,\cdot):=\Omega(\cdot,J \cdot)$ defines a 
pseudo-K\"ahler metric on $\mathfrak{g}$. By Section~\ref{sect: integrability}, the existence of a compatible pair on a Lie algebra is necessary 
for the existence of an integrable Born structure on the respective Lie algebra. 

\subsection{Dimension 4}
\label{sect:nil3 plus r}

According to~\cite[Section~9.4]{kuennethgeom}, the only non-Abelian $4$-di\-men\-sio\-nal Lie algebra admitting a K\"unneth 
structure is $\mathfrak{nil}_3\oplus\mathbb{R}$, the Lie algebra of $Nil^3\times \mathbb{R}$. 
This has a basis $\{e_1,e_2,e_3,e_4\}$ with the only non-zero bracket relation $[e_1,e_2]=e_3$.

This Lie algebra not only has K\"unneth structures, it even has 
a hypersymplectic structure, so that we can apply the construction of Subsection~\ref{ss:hyper} to it.
\begin{theorem}
There is an $S^1$-family of integrable Born structures on $\mathfrak{nil}_3\oplus \mathbb{R}$.
\end{theorem}
\begin{proof}
It is proved in \cite[Section~9.5]{kuennethgeom} that the symplectic forms
\begin{align*}
    \alpha =& \alpha_{14}-\alpha_{23}\\
    \beta=& -\alpha_{13}-\alpha_{24}\\
    \omega =& -\alpha_{13}+\alpha_{24}
\end{align*}
define a hypersymplectic structure
\begin{center}
        \begin{tikzcd}[column sep=small]& \omega\arrow[dl,"A"']\arrow{dr}{B} & \\\alpha \arrow{rr}{J} & & \beta
     \end{tikzcd}
    \end{center}
on $\mathfrak{nil}_3\oplus\mathbb{R}$. Here, the $\{\alpha_i\}$ denote the dual basis of $\{e_i\}$ and we use the abbreviation $\alpha_{ij}=\alpha_i\wedge\alpha_j$. The recursion operators of the hypersymplectic structure act on the basis $\{e_i\}$ by 
\begin{align*}
    Ae_1=e_2,\quad Ae_2=e_1,\quad Ae_3=-e_4,\quad Ae_4=-e_3\\
    Be_1=e_1,\quad Be_2=-e_2,\quad Be_3=e_3,\quad Be_4=-e_4,\\
    Je_1=e_2,\quad Je_2=-e_1,\quad Je_3=-e_4,\quad
    Je_4=-e_3.
\end{align*}
It follows that 
\begin{align*}
    B^*\alpha_1=\alpha_1,\quad B^*\alpha_2=-\alpha_2,\quad B^*\alpha_3=\alpha_3,\quad B^*\alpha_4=-\alpha_4.
\end{align*}
Therefore, the hypersymplectic metric $g(\cdot,\cdot)=\alpha(\cdot,B\cdot)$ can be expressed in terms of the dual basis $\{\alpha_i\}$ as
\begin{align*}
    g&=\lb\alpha_1\otimes B^*\alpha_4 - \alpha_4\otimes B^*\alpha_1\rb - \lb\alpha_2\otimes B^*\alpha_3 - \alpha_3\otimes B^*\alpha_2\rb\\ 
    &=-\lb\alpha_1\otimes\alpha_4+\alpha_4\otimes\alpha_1+ \alpha_2\otimes\alpha_3+\alpha_3\otimes\alpha_2\rb.
\end{align*}
Let $\tilde{J}$ be the endomorphism defined by
\begin{align*}
\tilde{J}e_1=e_2,\quad \tilde{J}e_2=-e_1,\quad \tilde{J}e_3=e_4,\quad \tilde{J}e_4=-e_3.    
\end{align*}
Then $\tilde{J}^2=-\id$ by construction and $\tilde{J}$ anti-commutes with $A$ and $B$. Since the only non-vanishing bracket relation of $\mathfrak{nil}_3\oplus\mathbb{R}$ is $[e_1,e_2]=e_3$ and $e_1$ and $e_2$ are mapped to each other under $\tilde{J}$, the only component of the Nijenhuis tensor of $\tilde{J}$ that does not vanish trivially is $N_{\tilde{J}}(e_1,e_2)$. For this, we observe
\begin{align*}
N_{\tilde{J}}(e_1,e_2)&=[\tilde{J}e_1,\tilde{J}e_2]-[e_1,e_2]-\tilde{J}\lb[\tilde{J}e_1,e_2]+[e_1,\tilde{J}e_2]\rb\\
&= -[e_2,e_1]-[e_1,e_2]-\tilde{J}\lb [e_2,e_2]-[e_1,e_1]\rb\\ &= 0.
\end{align*}
Moreover, we have
\begin{align*}
    \tilde{J}^*\alpha_1=-\alpha_2,\quad \tilde{J}^*\alpha_2=\alpha_1,\quad \tilde{J}^*\alpha_3=-\alpha_4,\quad \tilde{J}^*\alpha_4=\alpha_3
\end{align*}
and it follows that
\begin{align*}
    \tilde{J}^*g&=-\left\{ \tilde{J}^*\alpha_1\otimes \tilde{J}^*\alpha_4+\tilde{J}^*\alpha_4\otimes\tilde{J}^*\alpha_1+\tilde{J}^*\alpha_2\otimes\tilde{J}^*\alpha_3+\tilde{J}^*\alpha_3\otimes\tilde{J}^*\alpha_2\right\}\\
    &=-\left\{ -\alpha_2\otimes\alpha_3+\alpha_3\otimes(-\alpha_2)+\alpha_1\otimes(-\alpha_4)-\alpha_4\otimes\alpha_1\right\}\\
    &=-g.
\end{align*}
Therefore, by Theorem \ref{thm: s1-fam int bs} we obtain an $S^1$-family of integrable Born structures \begin{center}
    \begin{tikzcd}[column sep=small]& g\arrow[dl,"I_{\theta}"']\arrow{dr}{\tilde{B}_{\theta}} & \\\beta_{\theta} \arrow{rr}{-\tilde{J}} & & h_{\theta}
     \end{tikzcd},
\end{center}
where
\begin{align*}
    \beta_{\theta}=-\sin(\theta)\alpha + \cos(\theta)\beta, \quad  I_{\theta}=\cos(\theta)A + \sin(\theta)B
\end{align*}
and $\tilde{B}_{\theta}=\tilde{J}I_{\theta}$.
\end{proof}
Thus the Lie group $Nil^3\times \mathbb{R}$ carries left-invariant integrable Born structures that descend to all 
the associated nilmanifolds. These nilmanfolds are exactly the non-trivial principal $T^2$-bundles over $T^2$.

It was shown in~\cite[Theorem~10.27]{kuennethgeom} that the infra-nilmanifolds of $Nil^3\times \mathbb{R}$ also
carry K\"unneth structures, although they are not hypersymplectic since the hypersymplectic structure does not
descend. These infra-nilmanifolds cannot carry integrable Born structures since they are not complex manifolds.
 The genuine infra-nilmanifolds of type $Nil^3\times \mathbb{R}$ have even first Betti number $b_1=2$. 
 It is a fact that was first known from Kodaira's classification of compact complex surfaces and later proved 
 directly by Buchdahl~\cite{buchdahl2003compact} and Lamari~\cite{lamari1999courants} that compact 
 complex surfaces with even first Betti number are K\"ahler. However, since infra-nilmanifolds are finitely covered by nilmanifolds, such a K\"ahler structure would lift to these nilmanifolds, which is impossible, since they have first Betti number $3$ and so cannot be K\"ahler. 

It was shown in~\cite[Example~9.46]{kuennethgeom} that the Lie group $Sol^3\times \mathbb{R}$,
where $Sol^3$ is the solvable non-nilpotent Thurston geometry, admits left-invariant K\"unneth structures.
These cannot be upgraded to integrable Born structures because again the closed manifolds $\Gamma \setminus (Sol^3\times \mathbb{R})$ carrying 
such a geometry have even first Betti number, and so would be K\"ahler by the result mentioned above.
The lattice $\Gamma$ would then be a solvable K\"ahler group that is not virtually nilpotent. Such K\"ahler groups 
cannot exist be a result of Delzant~\cite{delzant2010invariant}.

\subsection{Dimension 6}

The $6$-dimensional nilpotent Lie algebras carrying K\"unneth structures were classified in~\cite{hamilton2019bi}, 
see also the summary in~\cite[Section~9.4.3]{kuennethgeom}. The upshot is that there are $15$ different non-Abelian Lie algebras
that do have this structure. Similarly, the $6$-dimensional nilpotent Lie algebras carrying pseudo-K\"ahler structures
were classified by Cordero, Fern\'andez and Ugarte~\cite{cordero2004pseudo}, who found $13$ non-Abelian cases.
These $13$ Lie algebras all have K\"unneth structures by the result of~\cite{hamilton2019bi}, and so they are candidates 
for having integrable Born structures. The two examples from~\cite{hamilton2019bi} which have K\"unneth structures but 
do not have pseudo-K\"ahler 
structures do not have to be considered.

We will show now that some of the pseudo-K\"ahler Lie algebras from~\cite{cordero2004pseudo} actually have integrable
Born structures, whereas at least one of them does not. We do not decide this question in all cases, leaving the completion
of the classification for future research.

\begin{example}\label{nonex}
	The Lie algebra denoted $\mathfrak{h}_{15}$ in~\cite{cordero2004pseudo} cannot have an integrable Born structure, since it was shown 
	by Andrada~\cite{andrada2008} that it does not admit a complex product structure; cf.~Remark~\ref{rem:complexprod}.
	\end{example}

\begin{example}
The Lie algebra denoted by $\mathfrak{h}_{8}$ in~\cite{cordero2004pseudo}  is the direct sum $\mathfrak{nil}_3\oplus\mathbb{R}^3$.
Since we have seen that $\mathfrak{nil}_3\oplus\mathbb{R}$ has integrable Born structures, the same is true in this case,
since we can just sum the structures with the standard Born structure on $\mathbb{R}^2$ from Example~\ref{ex: cn}.
\end{example}

For more complicated examples we use the following:
\begin{lemma}
\label{lem: Nijenhuis}
Let $J$ be an almost complex structure with Nijenhuis tensor $N_J$. If $N_J(X,Y)$ vanishes, so do $N_J(JX,JY)$, $N_J(JX,Y)$ and $N_J(X,JY)$.
\end{lemma}
\begin{proof}
This follows from
\begin{align*}
    N_J(JX,JY)=&[J^2X,J^2 Y]-[JX,JY]-J\lb [J^2X,JY]+[JX,J^2 Y]\rb\\
    =&-N_J(X,Y)\\
    N_J(JX,Y)=&[J^2 X,Y]-[JX,Y]-J\lb [J^2 X,Y]+[JX,JY]\rb\\
    =&-J\cdot N_J(X,Y).
\end{align*}
\end{proof}

\begin{example}
The Lie algebra denoted by $\mathfrak{h}_{4}$ in~\cite{cordero2004pseudo} is not a direct sum of Lie algebras of lower dimension.
It has a basis $e_1 , \ldots , e_6$ for which the only non-zero commutators are $[e_1 , e_2] = -e_5$ and $[e_1, e_4]=[e_2,e_3]=-e_6$.
We shall work with the dual $1$-forms $\alpha_i$. They are closed for $i\leq 4$ and satisfy $\dd\alpha_5=\alpha_1\wedge\alpha_2$
and $\dd\alpha_6=\alpha_1\wedge\alpha_4 + \alpha_2\wedge\alpha_3$. We will abbreviate $\alpha_{ij}=\alpha_i\wedge\alpha_j$ and 
$\alpha_{ijk}=\alpha_i\wedge\alpha_j\wedge\alpha_k$.

The $2$-form $\omega=\alpha_{13}+\alpha_{26}+\alpha_{45}$ on $\mathfrak{h}_4$ is non-degenerate and since
\begin{align*}
    \dd\omega=-\alpha_{214}-\alpha_{223}-\alpha_{412}=0
\end{align*}
it is also closed and therefore symplectic. 

The complementary subspaces
\begin{align*}
    \mathfrak{g}_+:=\langle e_1,e_2,e_5\rangle,\quad \mathfrak{g}_-:=\langle e_3,e_4,e_6\rangle
\end{align*}
of $\mathfrak{h}_4$ are Lagrangian for $\omega$. Moreover, since $[e_1,e_2]=-e_5$, they are also integrable. Therefore, $(\omega,\mathfrak{g}_+,\mathfrak{g}_-)$ is a K\"unneth structure.

The endomorphism $J$ defined by
\begin{align*}
    Je_1=-2 e_3,\quad
    J e_2=-e_4,\quad 
    Je_3= \frac{1}{2}e_1, \quad
    Je_4=e_2,\quad
    Je_5=e_6,\quad
    Je_6=-e_5.
\end{align*}
squares to $-Id$. To see that it is integrable, we compute
\begin{align*}
    N_J(e_1,e_2)=&[-2e_3,-e_4]-[e_1,e_2]-J\lb [-2e_3,e_2]+[e_1,-e_4]\rb\\
    =&e_5+Je_6=0 \ .
\end{align*}
By Lemma \ref{lem: Nijenhuis} this implies that also $N_J(e_1,e_4)=N_J(e_2,e_3)=N_J(e_3,e_4)=0$. Moreover, since $N_J(e_1,Je_1)=0=N_J(e_2,Je_2)$, we find $N_J(e_1,e_3)=0=N_J(e_2,e_4)$. Furthermore, since brackets with $e_5$ and $e_6$ vanish in $\mathfrak{h}_4$, the remaining components of the Nijenhuis tensor vanish trivially. 

We will now show that $J$ is compatible with the K\"unneth structure $(\omega,\mathfrak{g}_+,\mathfrak{g}_-)$. The subalgebras $\mathfrak{g}_{\pm}$ are clearly interchanged by $J$. It remains to show that $\omega$ is compatible with $J$. For this, we observe that $J$ acts on the dual basis as
\begin{align*}
    J^*\alpha_1=\frac{1}{2}\alpha_3,\quad J^*\alpha_2=\alpha_4,\quad J^*\alpha_3=-2\alpha_1,\quad J^*\alpha_4=-\alpha_2,\quad J^*\alpha_5=-\alpha_6,\quad J^*\alpha_6=\alpha_5.
\end{align*}
It follows that 
\begin{align*}
J^*\omega=-\alpha_{31}+\alpha_{45}+\alpha_{26}=\omega 
\end{align*}
and, hence, $(\omega,J)$ is a compatible pair. We conclude that $(\omega,\mathfrak{g}_+,\mathfrak{g}_-)$ is a K\"unneth structure compatible with the complex structure $J$ and by the discussion in Section~\ref{sect: kunneth and born} this defines an integrable Born structure. 
\end{example}

\begin{example}
	The Lie algebra $\mathfrak{h}_9$ in~\cite{cordero2004pseudo} is three step nilpotent and again, cannot be written as a direct sum of Lie algebras of lower dimension. It has a basis $e_1,\ldots,e_6$ for which the only non-zero commutators are $[e_1,e_2]=-e_4$, $[e_1,e_4]=-e_6$ and $[e_2,e_5]=-e_6$. As before, we are working with the dual 1-forms $\alpha_i$ which satisfy $\dd\alpha_i=0$ for $i\leq4$, $\dd\alpha_5=\alpha_{12}$ and $\dd\alpha_6=\alpha_{14}+\alpha_{25}$.
	
The non-degenerate 2-form $\omega=\alpha_{13}+4\alpha_{26}-4\alpha_{45}$ on $\mathfrak{h}_9$ is closed:
\begin{align*}
	\dd\omega=-4\alpha_{214}-4\alpha_{223}+4\alpha_{412}=0  
\end{align*}
and hence symplectic.

Moreover, the abelian subalgebras
\begin{align*}
		\mathfrak{g}_+=\langle e_1,e_5,e_6\rangle,\quad \mathfrak{g}_-=\langle e_2,e_3,e_4\rangle
\end{align*}	
of $\mathfrak{h}_9$ are Lagrangian for $\omega$ and therefore, $(\omega,\mathfrak{g}_+,\mathfrak{g}_-)$ is a K\"unneth structure.
	
Furthermore, we consider the endomorphism $J$ of $\mathfrak{h}_9$ defined by	
\begin{align*}
	Je_1=-e_2,\quad Je_2=e_1,\quad Je_3=-\frac{1}{4}e_6,\quad Je_4=-e_5,\quad Je_5=e_4, \quad Je_6=4e_3.
\end{align*}
Clearly, $J^2=-\id$. Moreover, we find
\begin{align*}
	N_J(e_1,e_3)=&[-e_2,-\frac{1}{4}e_6]-[e_1,e_3]-J\lb [-e_2,e_3]+[e_1,-\frac{1}{4} e_6]\rb=0\\
	N_J(e_1,e_4)=&[-e_2,-e_5]-[e_1,e_4]-J\lb [-e_2,e_4]+[e_1,-e_5]\rb\\
	=&-e_6+e_6=0\\
	N_J(e_3,e_4)=&[-\frac{1}{4}e_6,-e_5]-[e_3,e_4]-J\lb [e_3,-e_5]+[-\frac{1}{4}e_6,e_4]\rb=0.
\end{align*}
Using Lemma \ref{lem: Nijenhuis} it follows that the remaining components of the Nijenhuis tensor also vanish and, hence, $J$ is integrable. 
	
The complex structure $J$ clearly interchanges the subalgebras $\mathfrak{g}_{\pm}$. To see that $J$ it is compatible with $\omega$ we observe that
	\begin{align*}
		J^*\alpha_1=\alpha_2,\quad J^*\alpha_2=-\alpha_1,\quad J^*\alpha_3=4\alpha_6,\quad J^*\alpha_4=\alpha_5,\quad J^*\alpha_5=-\alpha_4,\quad J^*\alpha_6=-\frac{1}{4}\alpha_3.
	\end{align*}
	This yields
	\begin{align*}
		J^*\omega=4\alpha_{26}+\alpha_{13}+4\alpha_{54}=\omega.
	\end{align*}
	Hence, the K\"unneth structure $(\omega,\mathfrak{g}_+,\mathfrak{g}_-)$ is compatible with the complex structure $J$ and following the discussion in Section~\ref{sect: kunneth and born} this defines an integrable Born structure on $\mathfrak{h}_9$.
\end{example}

\bigskip

\bibliographystyle{amsplain}

\begin{thebibliography}{99}

\bibitem{andrada2008}
A. Andrada,
	{\em Complex product structures on $6$-dimensional nilpotent Lie algebras},
	Forum Math. {\bf 20} (2008), 285--315

\bibitem{AS}
A. Andrada and S. Salamon, {\em Complex product structures on Lie algebras}, Forum Math.~{\bf 17} (2005), 261--295

\bibitem{bande2008geometry}
G. Bande and D. Kotschick, {\em The geometry of recursion operators}, Comm.~math.~phys.~{\bf 280} (2008), 737--749


\bibitem{born1938suggestion}
M. Born, {\em A suggestion for unifying quantum theory and relativity}, Proc.~Royal Soc.~London. Series A. Mathematical and Physical Sciences {\bf 165} (1938), 291--303

\bibitem{buchdahl2003compact}
N. Buchdahl, {\em Compact K{\"a}hler surfaces with trivial canonical bundle}, Ann.~of Global Anal.~Geom.~{\bf 23} (2003), 189--204


\bibitem{cordero2004pseudo}
L. A. Cordero, M. Fern{\'a}ndez and L. Ugarte, 
{\em Pseudo-K{\"a}hler metrics on six-dimensional nilpotent Lie algebras}, J.~Geom.~Phys.~{\bf 50} (2004), 115--137


\bibitem{Cr} 
V.~Cruceanu, P.~Fortuny and P.~M.~Gadea, {\em A survey on paracomplex geometry}, Rocky Mountain J.~Math.~{\bf 26} (1996), 83--115.

\bibitem{delzant2010invariant}
T. Delzant, {\em L'invariant de Bieri--Neumann--Strebel des groupes fondamentaux des vari{\'e}t{\'e}s k{\"a}hl{\'e}riennes}, Math.~Ann.~{\bf 348} (2010), 119--125

\bibitem{GS} 
F.~Etayo Gordejuela and R.~Santamar{\'{\i}}a, {\em The canonical connection of a bi-Lagrangian manifold}, 
Journ.~of Phys.~A {\bf 34} (2001), 981--987.

\bibitem{freidel2014born}
L. Freidel, R. G. Leigh and D. Minic, {\em Born reciprocity in string theory and the nature of spacetime}, Physics Letters B {\bf 730} (2014), 302--306

\bibitem{freidel2014quantum}
L. Freidel, R.  Leigh, and D. Minic, {\em Quantum Gravity, Dynamical Phase Space and String Theory}, International Journal of Modern Physics D {\bf 23} (2014) . 

\bibitem{freidel2015metastring}
L. Freidel, R. Leigh, D. Minic, {\em Metastring Theory and Modular Space-time}, Journal of High Energy Physics (2015).  

\bibitem{freidel2017generalised}
L. Freidel, F. J. Rudolph and D. Svoboda, {\em Generalised kinematics for double field theory}, Journal of High Energy Physics {\bf 2017} (2017), 1--35

\bibitem{freidel2019unique}
L. Freidel, F. J. Rudolph and D. Svoboda, {\em A unique connection for Born geometry}, Communications in Mathematical Physics {\bf 372} (2019), 119--150


\bibitem{hamilton2019bi}
M. J. D. Hamilton, {\em Bi-Lagrangian structures on nilmanifolds}, Journal of Geometry and Physics {\bf 140} (2019), 10--25

\bibitem{kuennethgeom}
M. J. D. Hamilton and D. Kotschick, {\sl K\"unneth geometry}, Cambridge Univ.~Press 2024.

\bibitem{H1} 
H.~Hess, {\em Connections on symplectic manifolds and geometric quantization}, Lect. Notes in Math. {\bf 836} (1980), 153--166.

\bibitem{H2} 
H.~Hess, {\em On a geometric quantization scheme generalizing those of Kostant-Souriau and Czyz}, Lect.~Notes in Physics {\bf 139} (1981), 1--35.


\bibitem{hitchin1990hypersymplectic}
N. Hitchin, {\em Hypersymplectic quotients}, Acta Acad. Sci. Tauriensis {\bf 124} (1990), 169--180


\bibitem{hu2019zro}
S.~Hu, R.~Moraru and D.~Svoboda,
{\em Commuting Pairs, Generalized para-K\"ahler Geometry and Born Geometry}, arXiv preprint
[arXiv:1909.04646 [hep-th]] (2019).


\bibitem{ivanov2005parahermitian}
S. Ivanov, and S. Zamkovoy, {\em Parahermitian and paraquaternionic manifolds}, Diff.~Geom.~Appl.~{\bf 23} (2005), 205--234.

\bibitem{lamari1999courants}
A. Lamari, {\em Courants k{\"a}hl{\'e}riens et surfaces compactes}, Ann.~de l'Inst.~Fourier {\bf 49} (1999), 263--285


\bibitem{libermann1952structures}
P. Libermann, {\em Sur les structures presque quaternioniennes de deuxi\`eme esp\`ece}, C. R. Acad. Sci. Paris {\bf 234} (1952), 1030--1032.

\bibitem{marotta2018parahermitian}
V. Marotta, R. Szabo, {\em Para-Hermitian Geometry, Dualities and Generalized Flux Backgrounds}, Fortschritte der Physik {\bf 67} (2018).

\bibitem{marotta2021born}
V. E. Marotta, R. J. Szabo, {\em Born sigma-models for para-{H}ermitian manifolds and generalized {T}-duality}, Reviews in Mathematical Physics {\bf 33} (2021)


\bibitem{svoboda2020born}
D. Svoboda, {\em Born Geometry}, PhD thesis (2020)



\bibitem{vukmirovic2003quaternionic}
S. Vukmirovic, {\em Para-quaternionic reduction}, arXiv preprint math/0304424, 2003.

\bibitem{zamkovoy2006geometry}
S. Zamkovoy, {\em Geometry of paraquaternionic K{\"a}hler manifolds with torsion}, Journal of Geometry and Physics {\bf 57} (2006), 69--87.


\end{thebibliography}

\end{document}